\documentstyle[epsf,psfig,amsmath,12pt]{article}
\input epsf
\newcommand{\rz}{\ifmmode{{\rm I} \hskip -2pt {\rm R}}
    \else{\hbox{$I\hskip -2pt R$}}\fi} %
\newcommand{\dz}{\ifmmode{{\rm I} \hskip -2pt {\rm D}}
    \else{\hbox{$I\hskip -2pt R$}}\fi}

\newcommand{\bu}{{\bf u}} 
\newcommand{\obu}{\overline{\bf u}} 
\newcommand{\oS}{\overline S}

\newcommand{\bff}{{\bf f}} 

\newcommand{\bx}{{\bf x}}  
\newcommand{\bk}{{\bf k}}  
\def\PP{{{\rm l}\kern - .15em {\rm P} }}
\def\PN2{{\PP_{N}-\PP_{N-2}}}

\begin{document}

\title{Large Eddy Simulation of Turbulent Channel Flows by the Rational LES Model}
\author{T. Iliescu and P. Fischer}

\maketitle

\begin{abstract}  The rational large eddy simulation (RLES) model is applied
to turbulent channel flows.   This approximate deconvolution 
model is based on a rational (subdiagonal 
Pad\'e) approximation of the Fourier transform of the Gaussian filter
and is proposed as an alternative to the gradient 
(also known as the nonlinear or tensor-diffusivity) model.
We used a spectral element code to perform large eddy simulations of
incompressible channel flows at
Reynolds numbers based on the friction velocity and the channel half-width 
$Re_{\tau}=180$ and $Re_{\tau}=395$.
We compared the RLES model with the gradient model.
The RLES results showed a clear improvement over those corresponding to the
gradient model, comparing well with the fine direct numerical simulation.
For comparison, we also present results corresponding to a classical subgrid-scale 
eddy-viscosity model such as the standard Smagorinsky model.
\end{abstract}

\newpage

\section{Introduction}
Large eddy simulation (LES) is one of the most successful techniques in the
numerical simulation of turbulent flows.
Contrary to the direct numerical simulation (DNS), which tries to capture
all the scales in the flow, LES aims at resolving only the large-scale flow
features.
The large scales are defined by means of a filtering operation: the 
Navier-Stokes equations are convolved with a spatial filter, and the 
resulting filtered  variables become the variables of interest in LES.
Thus, a good LES model should be able to compute an accurate approximation
of the filtered variables.

An essential challenge in LES is the modeling of the subgrid-scale (SGS) stresses,
representing the interactions between the large (above the filter width)
and small (below the filter width) scales in the filtered Navier-Stokes
equations.
A remarkable research effort has led to a wide variety of SGS models,
surveyed, for example, in \cite{MK00}, \cite{LM96}, and \cite{Sag01}.

Arguably the most popular class of LES models is the eddy-viscosity type,
based on (variants of) the Smagorinsky model \cite{Sma63}.
The main feature of the eddy-viscosity models is that they properly 
transfer kinetic energy (by inviscid processes) from large scales to 
smaller and smaller scales, until this energy is 
dissipated  through viscous effects.
These models have several limitations, however, including poor correlation
coefficients in a priori tests \cite{CFR79}, \cite{BFR80} and inability 
to provide backscatter.
Some of these limitations are circumvented by using a dynamic procedure in
calculating the Smagorinsky constant, yielding the dynamic subgrid-scale 
eddy-viscosity model introduced by Germano et al. \cite{GPMC91},
and used in many studies \cite{Pio93}, \cite{Hor97}.

Another class of SGS models is the scale-similarity one.
The scale-similarity model, introduced by Bardina et al. \cite{BFR80},
postulates that the full structure of the velocity field at scales below the
filter width is similar to that at scales above the filter width.
A priori tests \cite{BFR80} show high correlations between real and modeled 
stresses.  
Another realistic feature of the scale-similarity model is that it produces 
backscatter.
In a posteriori tests, however, the scale-similarity model does not dissipate 
enough energy and typically leads to inaccurate results.
As a remedy, Bardina et al. \cite{BFR80} added a dissipative Smagorinsky 
term.  
The resulting model, known as the mixed model, combines the strengths of
both the scale-similarity and the Smagorinsky model.
The dynamic procedure has been successfully applied to both the pure and the 
mixed scale-similarity model, yielding improved results \cite{SPB99}.

A different class of SGS models consists of those models aimed at computing an 
improved SGS stress approximation by replacing the unknown unfiltered variables 
with approximately deconvolved filtered variables.
An inverse filtered model was first proposed by Shah and Ferziger 
\cite{SF95}.  This idea was formalized by Geurts \cite{Geu97} for the top hat 
filter.
Kuerten et al. \cite{KGVG99} used the approximate inverse to improve the 
computable estimates in the dynamic Smagorinsky model.
Another model in this class is the velocity estimation model of Domaradzki 
and Saiki \cite{DS97}, \cite{DL99}, \cite{LD99}.
Stolz and Adams \cite{SA99} developed the approximate deconvolution model, 
based on repeated application of the filter to approximately deconvolve 
the dependent variables \cite{SAK01a}, \cite{SAK01b}.

One popular model in this class is the gradient model (also known as the 
nonlinear or tensor-diffusivity model), which uses {\it explicit} filtering.
In addition to the {\it implicit} filtering due to the effective 
truncations (grid and numerical method), this LES model also assumes a regular 
explicit filter of prescribed shape and effective width larger than the grid
spacing.

The gradient model is based on a Taylor series approximation of the Fourier transform of 
the filter and aims at reconstructing the filtered-scale stress
due to explicit filtering.
The gradient model was developed in several steps.
First, in 1974 Leonard \cite{Leo74} proposed a model for the ``resolved scales'' 
${\overline{{\overline{\bu}} \ {\overline{\bu}}}}$ in the Reynolds stress tensor.
Next, in 1979 Clark, Ferziger, and Reynolds \cite{CFR79} used the same approach
to model the ``cross-terms'' 
${\overline{{\overline{\bu}}\bu'}} + {\overline{\bu' {\overline{\bu}}}}$.

The gradient model was tested a priori against experimental data (two-dimensional 
cuts) by Liu et al. \cite{LMK94}.
Borue and Orszag \cite{BO98} presented a detailed a priori analysis of the 
gradient model based on Gaussian-filtered DNS of homogeneous, isotropic decaying 
turbulence.
Also, Winckelmans et al. \cite{WWVJ01} presented several a priori tests for 
the gradient model and its dynamic version, again in the context of homogeneous,
isotropic decaying turbulence.
Similar tests have been performed by Carati et al. \cite{CWJ01}.
All the above a priori tests have shown high correlations.

In {\it a posteriori} tests, however, it was found that the gradient model 
does not dissipate enough energy.
Simulations with the pure gradient model appear to be unstable \cite{Vre95}.
Also, Liu, Meneveau, and Katz \cite{LMK94} reported problems near the wall,
where the pure gradient model's Reynolds stresses do not follow the $x_2^3$ 
behavior.
To stabilize the gradient model, Clark, Ferziger, and Reynolds 
\cite{CFR79} combined it with a Smagorinsky term, but the resulting mixed model
inherited the excessive dissipation of the Smagorinsky model.
A different approach was proposed by Liu et al. \cite{LMK94}, who supplied the 
gradient model with a ``limiter'' to prevent energy backscatter; 
this clipping procedure ensures that the model dissipates energy from large to
small scales.  This approach was also used in \cite{CW97}, \cite{CV98}.

From this point of view, the gradient model is similar to the scale-similarity
model: it shows high correlations in a priori tests, but it does not dissipate 
 enough energy in actual LES simulations: hence the need for extra viscosity
type terms (mixed models).
We note that, for both types of model, the best results in actual
LES simulations were obtained by using the dynamic mixed procedure \cite{VGK97},
\cite{WWVJ01}.
In fact, it has been noted before \cite{VGK97}, \cite{Hor97}, \cite{CWJ01}, 
\cite{{WWVJ01}} that there are strong ties between the gradient model and the
scale-similarity model: the first term in the Taylor series expansion of the
scale-similarity model is indeed the gradient model.
As noted by Winckelmans et al. \cite{WWVJ01}, however, the other terms in the 
expansion are different.
Thus, the gradient model is not identical to the scale-similarity model.

The model presented in this paper was introduced by Galdi and Layton \cite{GL00}
as an alternative to the gradient model.
They observed that the Taylor series approximation of the Fourier transform of the
Gaussian filter used in the derivation of the gradient model actually {\it increases}
the high wave number components, instead of damping them.
As an alternative to the Taylor series approximation, Galdi and Layton proposed a 
rational ((0,1) Pad\'e) approximation.
This rational approximation is consistent with the original approximated function
(which is a negative exponential): it attenuates the high wave number components.

In this paper, the resulting LES model, called in the sequel the rational LES 
(RLES) model, is applied to the numerical simulations of incompressible channel
flows at $Re_{\tau}=180$ and $Re_{\tau}=395$.

\section{The Rational LES Model}

The usual LES starts by convolving the Navier-Stokes equations (NSEs) with a 
spatial filter $g_{\delta}$.
Assuming that differentiation and convolution commute (which is true for 
homogeneous filters), the filtered NSEs read as follows:
\begin{eqnarray}\label{Filtered_NSEs}
{\overline{\bu}}_t + \nabla \cdot ({\overline{\bu \bu}}) 
- Re^{-1} \Delta{\overline{\bu}} + \nabla {\overline{p}} =  {\overline{\bff}},
\end{eqnarray}

\noindent where $\delta$ is the filter width and 
${\overline{\bu}}=g_{\delta}*{\bu}$ is the variable of interest.  
The filtered NSEs (\ref{Filtered_NSEs}) do not form a closed system, 
and a considerable research effort in LES 
research has been directed at modeling the stress
\begin{eqnarray}\label{stress} 
\tau = {\overline{\bu \bu}} - {\overline{\bu}} \ {\overline{\bu}}.
\end{eqnarray}

As mentioned by Carati et al. \cite{CWJ01}, this stress consists of 
a filtered-scale stress tensor, mainly due to filtering, and a 
subgrid-scale (SGS) stress tensor, mainly due to discretization.
One way of approximating the filtered-scale stress tensor 
is by using a Taylor series expansion in the wave number space to represent the
unknown full velocity in terms of the filtered velocity.
This approach was first used by Leonard \cite{Leo74}, and it was later espoused 
by Clark, Ferziger, and Reynolds \cite{CFR79}.
The resulting model, called the gradient, nonlinear, or tensor-diffusivity model,
was used in numerous studies \cite{Leo74}, \cite{CFR79}, \cite{CWJ01}, 
\cite{BO98}, \cite{WWVJ01}, \cite{VGK96}, \cite{LMK94}, \cite{Ald90},
\cite{KSF00}.

The gradient model is derived by using a Taylor series approximation to the
Fourier transform of the Gaussian filter
\begin{eqnarray}\label{taylor.approx} 
{\widehat{g_{\delta}}}(\bk) 
= e^{-\frac{\delta^2 |\bk|^2}{4 \gamma}}
\approx 1 - \frac{\delta^2 |\bk|^2}{4 \gamma} +O(\delta^4),
\end{eqnarray}

\noindent and for its inverse
\begin{eqnarray}\label{taylor.approx.inverse} 
\frac{1}{{\widehat{g_{\delta}}}(\bk)} 
= e^{\frac{\delta^2 |\bk|^2}{4 \gamma}}
\approx 1 + \frac{\delta^2 |\bk|^2}{4 \gamma} +O(\delta^4).
\end{eqnarray}

\noindent Decomposing $\bu$ into its average and its turbulent fluctuations
\begin{eqnarray}\label{decomposition.u} 
\bu = \obu + \bu',
\end{eqnarray}

\noindent and taking first the average and then the Fourier transform of the
above relation, we get
\begin{eqnarray}\label{uprime1} 
{\widehat{\bu'}} = \left( \frac{1}{{\widehat{g_{\delta}}}} - 1 \right) {\widehat{\obu}},
\end{eqnarray}

\noindent and thus,
\begin{eqnarray}\label{uprime2}
{\widehat{\bu}} = \frac{1}{{\widehat{g_{\delta}}}} {\widehat{\obu}},
\end{eqnarray}
 
\noindent where ${\widehat \bu}$ denotes the Fourier transform of $\bu$.

\noindent By taking the inverse Fourier transform and using 
(\ref{taylor.approx.inverse}), we get
\begin{eqnarray}
\bu \approx \obu + \frac{\delta^2}{4 \gamma} \Delta \obu.
\end{eqnarray}
 
\noindent By plugging the above into (\ref{stress}), using 
(\ref{taylor.approx}) and the same technique as above, 
simplifying, and dropping out the terms of $O(\delta^4)$,
we get the gradient model
\begin{eqnarray}\label{gradient.model}
\tau  = {\overline{\bu \bu}} - {\overline{\bu}} \ {\overline{\bu}}
\approx \frac{\delta^2}{2 \gamma} \nabla \obu \nabla \obu,
\end{eqnarray}
 
\noindent where 
\begin{eqnarray}  
(\nabla \obu \nabla \obu)_{i,j} = \sum_{l=1}^{d} \frac{\partial \obu_i}{\partial \bx_l} 
\frac{\partial \obu_j}{\partial \bx_l}.  
\end{eqnarray}

\begin{figure}
\centerline{
\hbox{
\psfig{figure=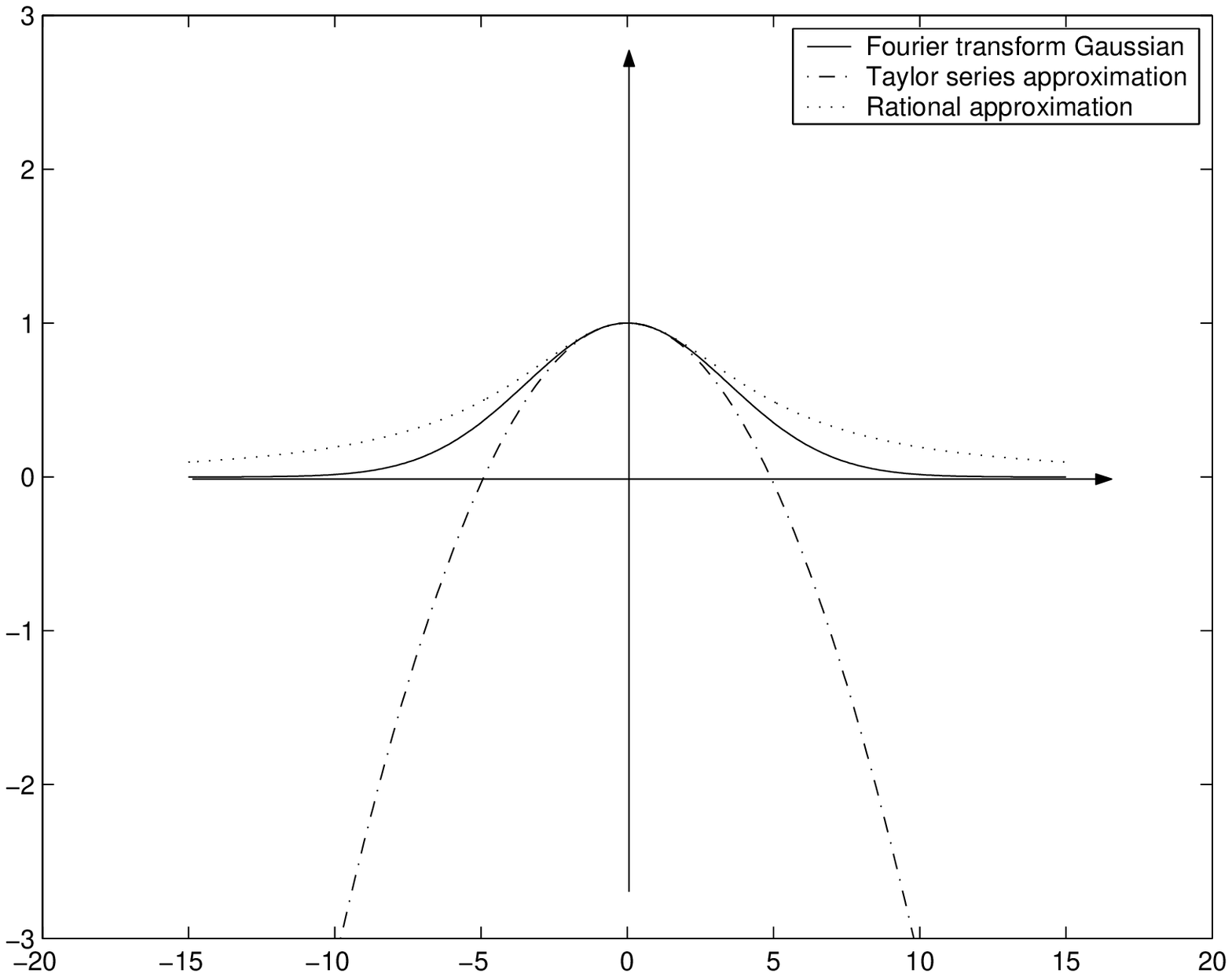,height=250pt,width=250pt,angle=0}}
}
\caption{\label{Rational.vs.Taylor} Approximations to the Fourier 
transform of the Gaussian filter: Rational (Pad\'e) vs. Taylor.}
\end{figure}

Noticing that the approximation by Taylor series of ${\widehat{g_{\delta}}}$ actually 
{\em increases} the high wave number components (see Figure \ref{Rational.vs.Taylor}), 
Galdi and Layton \cite{GL00} 
developed a new LES model based on a rational ((0,1) Pad\'e) approximation of 
${\widehat{g_{\delta}}}$, which preserves the decay of the high wave number components:
\begin{eqnarray}\label{rational.approx} 
{\widehat{g_{\delta}}}(\bk) 
= e^{-\frac{\delta^2 |\bk|^2}{4 \gamma}}
\approx \frac{1}{1 + \frac{\delta^2 |\bk|^2}{4 \gamma}} +O(\delta^4).
\end{eqnarray}

\noindent The resulting LES model, called the rational 
LES (RLES) model, reads as follows:
\begin{eqnarray}\label{rational.LES}  
\tau = \left[ \left( -\frac{\delta^2}{4 \gamma} \Delta + I \right)^{-1} 
\left( \frac{\delta^2}{2 \gamma} \nabla \obu \nabla \obu \right) \right].  
\end{eqnarray}

The inverse operator in (\ref{rational.LES}) acts as a smoothing operator
and represents the approximation of the convolution by the Gaussian filter
in the stress tensor $\tau$ in (\ref{stress}).

We note that differential filters have been proposed by 
Germano in \cite{Ger86}: Actually, one can think of (\ref{rational.LES})
as the stress tensor obtained by applying such a differential filter.
Mullen and Fischer used similar filters in \cite{MF99}.
Also, Domaradzki and Holm considered the Navier-Stokes-alpha
model (which contains an inverse operator similar to the one in 
(\ref{rational.LES})), in an LES framework \cite{DH00}. 

The mathematical analysis associated with the RLES model (\ref{rational.LES})
was presented in \cite{BGIL01}; the smoothing property of the inverse
operator in (\ref{rational.LES}) eliminated the necessity for using 
additional regularization operators (of eddy-viscosity type), as for 
the gradient model \cite{Col97}.
The first steps in the numerical analysis and validation of the 
RLES model (\ref{rational.LES}) were made in \cite{IJL01} and
\cite{IJLMT01}, respectively.

This paper presents numerical results for the RLES model (\ref{rational.LES}) 
applied to the 3D channel flow test problem at Reynolds numbers based on the wall 
shear velocity $Re_{\tau}=180$ and $Re_{\tau}=395$. 
Some preliminary work started in \cite{FI01}; it was significantly updated and
improved in the present paper.

\section{Numerical Setting}
\begin{figure}
\centerline{
\hbox{
\psfig{figure=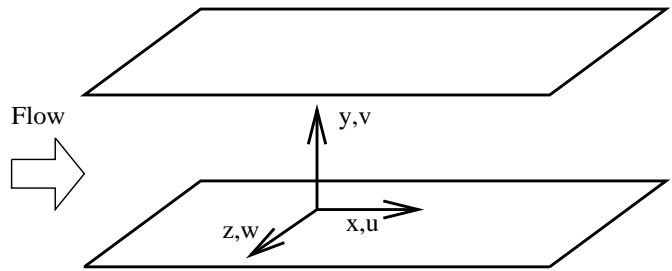,height=100pt,width=250pt,angle=270}}
}
\caption{\label{Channel} Problem setup for the channel flow.}
\end{figure}

The 3D channel flow (Figure \ref{Channel}) is one of the most popular test problems 
for the investigation of wall bounded turbulent flows \cite{MK82}, \cite{KMM87}.  
We used the fine DNS of Moser, Kim, and Mansour \cite{MKM99} as benchmark for 
our LES simulations.

We compared the RLES model (\ref{rational.LES}) with\\ 
\hspace*{0.5in} (I) the gradient model (\ref{gradient.model})
$\tau = \frac{\delta^2}{2 \gamma} \nabla \obu \nabla \obu$;\\ 
\hspace*{0.5in} (II) the Smagorinsky model $\tau = - (C_s \delta )^2 \ |\oS| \ \oS$;\\
\hspace*{0.5in} (III) a coarse DNS (no LES model),\\
\noindent where $\oS:=\frac{1}{2}(\nabla \obu + \nabla \obu^T)$ is the deformation tensor of 
the filtered field, $C_s = 0.1$ is the Smagorinsky constant, and $\gamma=6$ is the parameter
in the definition of the Gaussian filter.

The computational domain is periodic in the streamwise ($x$) and spanwise ($z$) directions, 
and the pressure gradient that drives the flow is adjusted dynamically to maintain 
a constant mass flux through the channel.  
The parameters used in the numerical simulations are given in Table \ref{table.param}
for the two Reynolds numbers considered ($Re_{\tau}=180$ and $Re_{\tau}=395$).
The filter width $\delta$ is computed as $\delta=\sqrt[3]{\Delta_x \ \Delta_z \ \Delta_y(y)}$, 
where $\Delta_x$ and $\Delta_z$ are the largest spaces between the 
Gauss-Lobatto-Legendre (GLL) points in the $x$ and $z$ directions, respectively, 
and $\Delta_y(y)$ is inhomogeneous and is computed as an interpolation function 
that is zero at the wall and is twice the normal mesh size for the elements 
in the center of the channel.

\begin{figure}
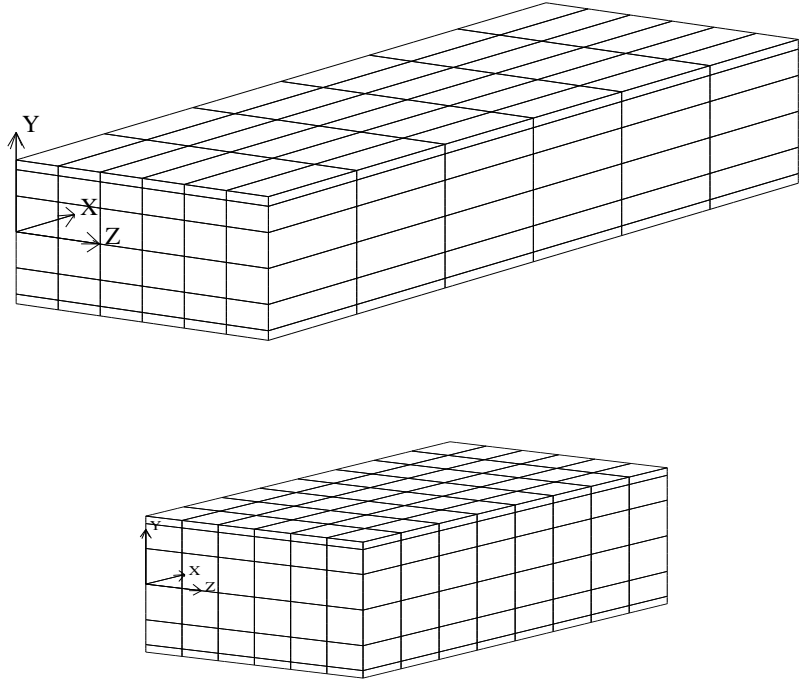

\centerline{
\hbox{
\psfig{figure=mesh180.epsi,height=130pt,width=300pt,angle=270}}
}\vspace*{0.5in}
\centerline{
\hbox{
\psfig{figure=mesh395.epsi,height=90pt,width=200pt,angle=270}}
}
\caption{\label{Mesh} Spectral element meshes: $Re_\tau=180$ (top),
                      and $Re_\tau=395$ (bottom).}
\end{figure}

The numerical simulations were performed by using a spectral element
code based on the $\PN2$ velocity and pressure spaces introduced
by Maday and Patera \cite{MP89}.
The domain is decomposed into spectral elements, 
as shown in Figure \ref{Mesh}.
Mesh spacing in the wall-normal direction ($y$) was chosen to
be roughly equivalent to a Chebychev distribution having the 
same number of points.
The velocity is continuous across element interfaces and is 
represented by $N$th-order tensor-product Lagrange polynomials based 
on the GLL points.  The pressure is discontinuous and is represented
by tensor-product polynomials of degree $N-2$.
Time stepping is based on an operator-splitting of the discrete
system, which leads to separate convective, viscous, and pressure
subproblems without the need for ad hoc pressure boundary conditions.
A filter, which removes 2\%--5\% of the highest velocity mode, is used to
stabilize the Galerkin formulation \cite{FM01};
the filter does not compromise the spectral accuracy.
Details of the discretization and solution algorithm are given in
\cite{Fis97}, \cite{FMT00}.

The initial conditions for the $Re_{\tau}=180$ simulations were obtained by 
superimposing a 2D Tollmien-Schlichting (TS) mode of 2\% amplitude and 
a 3D TS mode of 1\% amplitude on a parabolic mean flow (Poiseuille flow)
and integrating the flow for a long time (approximately 200 H/$u_{\tau}$,
where H is the channel's halfwidth and $u_{\tau}$ is the wall shear velocity)
on a finer mesh ($72 \times 73 \times 72$ mesh points).
The final field file was further integrated on the actual coarse LES
mesh ($36 \times 37 \times 36$ mesh points) for approximately 50 H/$u_{\tau}$
to obtain the initial condition for {\it all} four $Re_{\tau}=180$
simulations.  

The initial condition for the $Re_{\tau}=395$ case was obtained in a 
similar manner: We started with a field file corresponding to a 
$Re_{\tau}=180$ simulation, and we integrated it on a finer mesh
($96 \times 73 \times 72$ mesh points) for a long time
(approximately 50 H/$u_{\tau}$).
Then, we integrated the resulting flow on the actual coarser LES
mesh ($72 \times 55 \times 54$ mesh points) for another 
40 H/$u_{\tau}$, and the final field file was used as initial condition
for {\it all} four simulations.

\begin{figure}
\centerline{
\hbox{
\psfig{figure=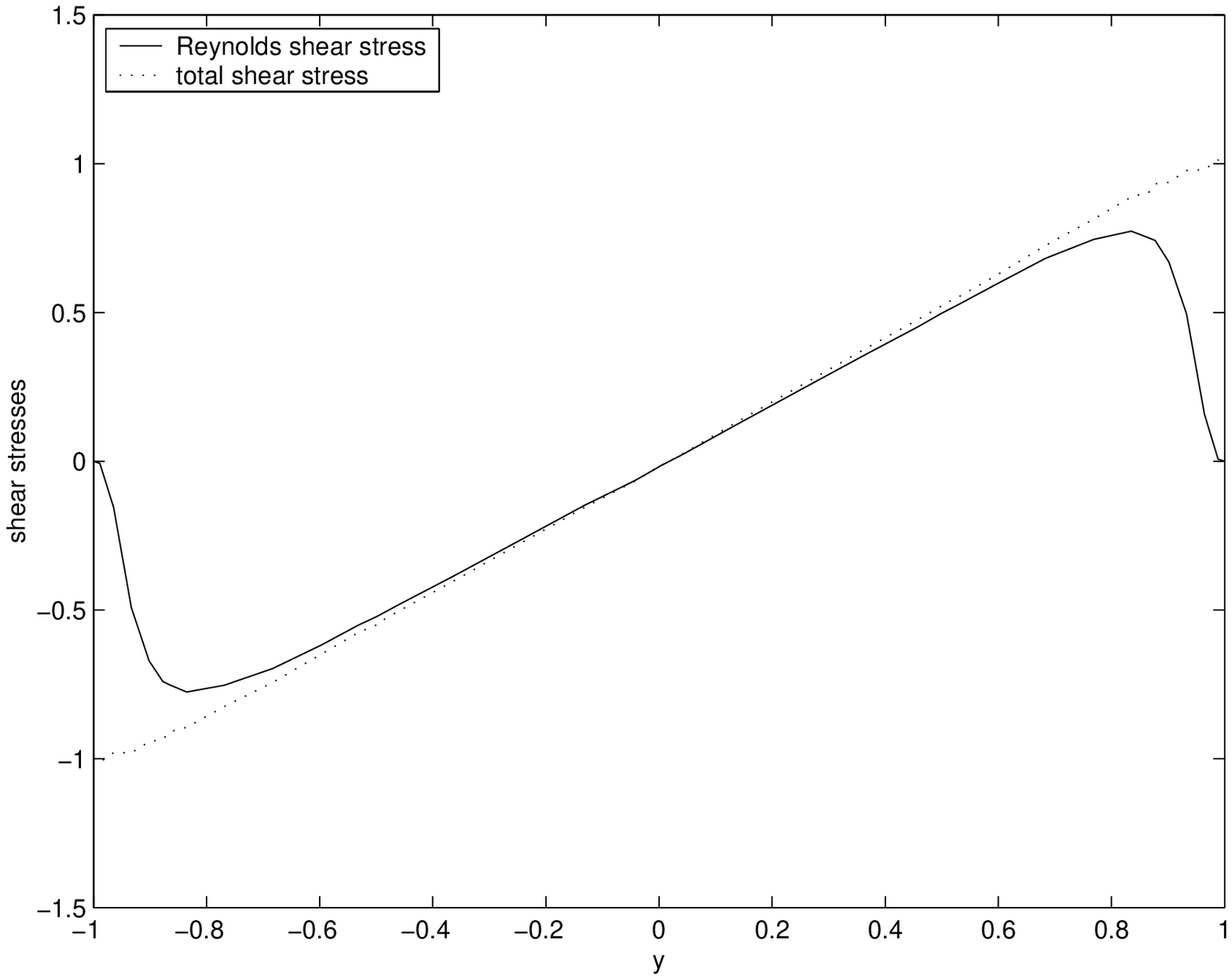,height=300pt,width=400pt,angle=0}}
}
\caption{\label{Figure0a} $Re_{\tau}=180$, linear total shear stress profile,
an indication that the statistically steady state was reached. 
}
\end{figure}

\begin{figure}
\centerline{
\hbox{
\psfig{figure=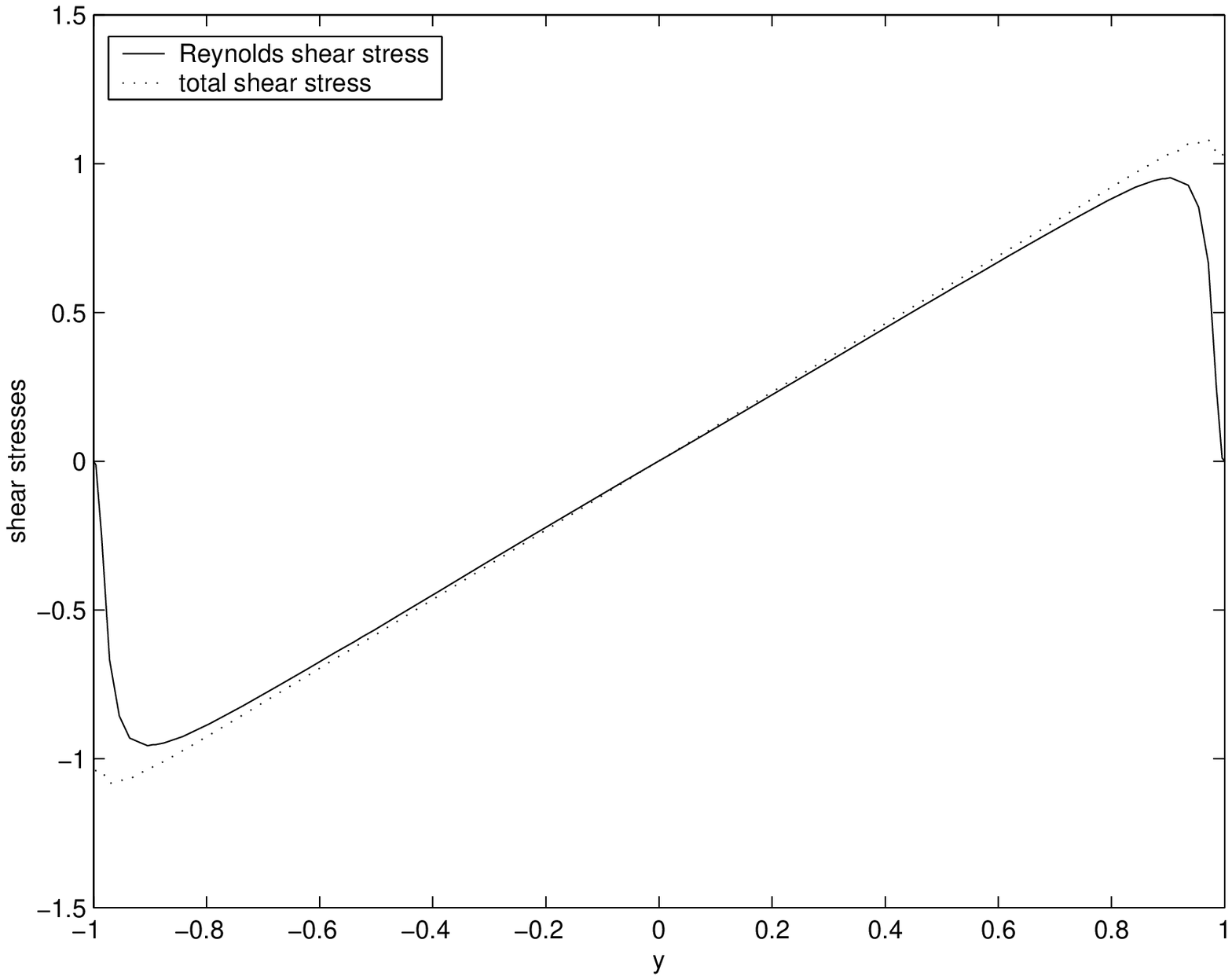,height=300pt,width=400pt,angle=0}}
}
\caption{\label{Figure0b} $Re_{\tau}=395$, linear total shear stress profile,
an indication that the statistically steady state was reached. 
}
\end{figure}

For each of the four simulations and for both $Re_{\tau}=180$ and 
$Re_{\tau}=395$, the flow was integrated further in time until the 
statistically steady state was reached (for approximately 
15 H/$u_{\tau}$).
The statistically steady state was identified by a linear
total shear stress profile (see Figure \ref{Figure0a} and 
Figure  \ref{Figure0b}).

The statistics were then collected over another 5 H/$u_{\tau}$
and contained samples taken after each time step
($\Delta t = 0.0002$ for $Re_{\tau}=180$ and 
$\Delta t = 0.00025$ for $Re_{\tau}=395$).
We also averaged over the two halves of the channel.

In our numerical experiments, we considered, as a first step, 
homogeneous boundary conditions for all LES models tested.

The numerical results include plots of the following 
time- and plane-averaged 
quantities normalized by the {\it computed} $u_{\tau}$: 
the mean streamwise velocity $\ll {\overline{u}}\gg / u_{\tau}$; 
the $x,y$-component of the Reynolds stress $\ll u' v'\gg / u_{\tau}^2$; 
and the rms values of the streamwise $\ll u' u'\gg / u_{\tau}^2$; 
wall-normal $\ll v' v'\gg / u_{\tau}^2$, and spanwise $\ll w' w'\gg / u_{\tau}^2$ 
velocity fluctuations, 
where $\ll \cdot\gg $ denotes 
time and plane ($xz$) averaging, the fluctuating quantities $f'$ are 
calculated as $f' = f - \ll f \gg $, and a ``$^{+}$'' superscript 
denotes the variable in wall units.

Note that in our simulations the bulk velocity $U_m$ was fixed
to match the corresponding one in \cite{MKM99} 
(see Table \ref{table.utau}), 
and the friction velocity $u_\tau$ was a result of the simulations.
Table 2 presents the {\it actual} values of $Re_\tau$ corresponding 
to the friction velocity $u_\tau$ computed for all four numerical
tests and two nominal Reynolds numbers.
We note that the friction velocity $u_{\tau}$ is within 1\%--2\%
of the nominal value, and, as a result, so is the actual $Re_{\tau}$.

\section{{\it A Posteriori} Tests for $Re_{\tau}=180$}
We ran {\em a posteriori} tests for 
the RLES model (\ref{rational.LES}), 
the gradient model (\ref{gradient.model}), 
the Smagorinsky model, and a coarse DNS (no model). 
We compared the corresponding results with the fine DNS simulation of 
Moser, Kim, and Mansour \cite{MKM99}.  

\begin{figure}
\centerline{
\hbox{
\psfig{figure=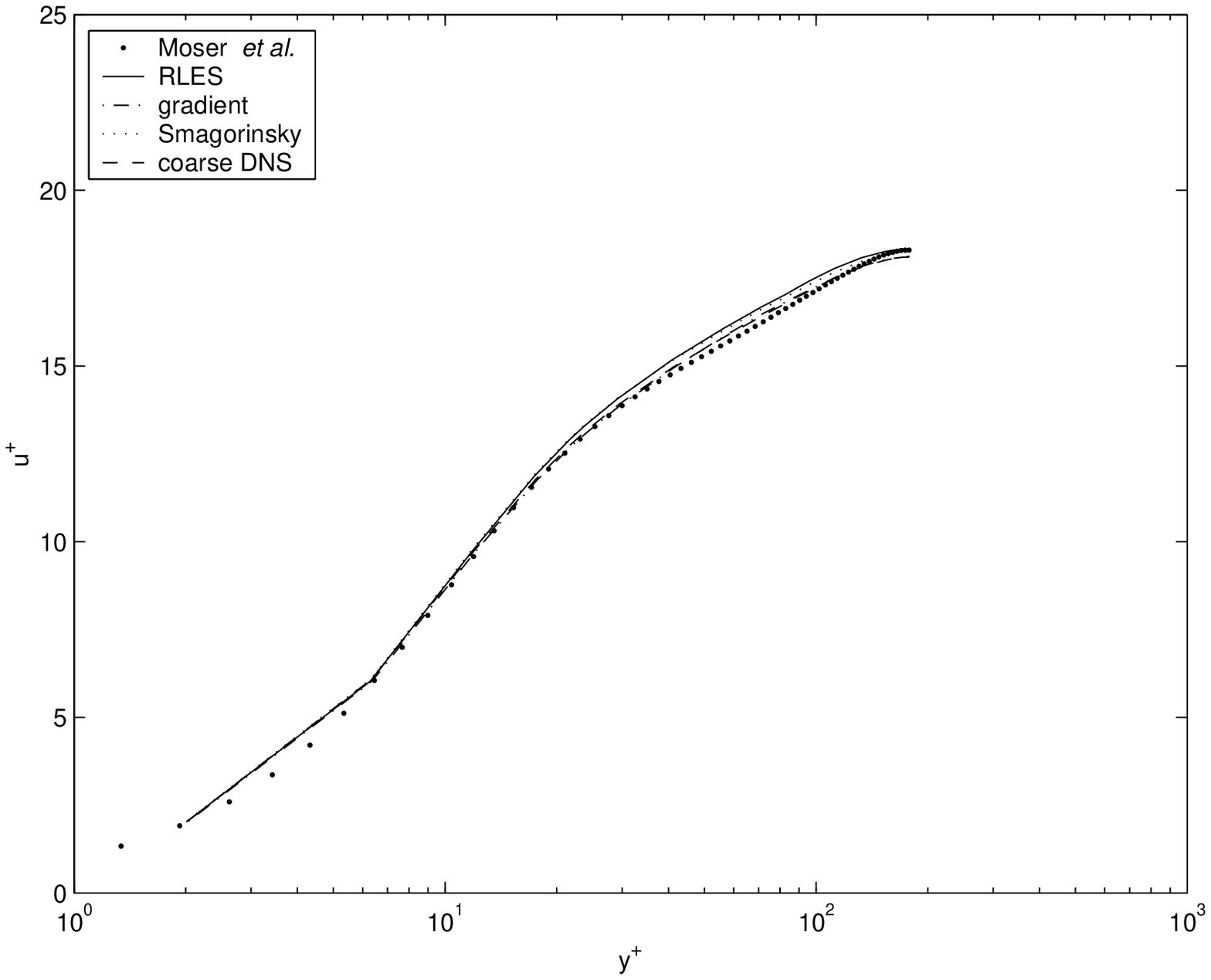,height=300pt,width=400pt,angle=0}}
}
\caption{\label{Figure1a} Mean streamwise velocity, $Re_{\tau}=180$.
We compared the RLES model (\ref{rational.LES}), the gradient model
(\ref{gradient.model}), the Smagorinsky model, and a coarse DNS,
with the fine DNS of Moser, Kim, and Mansour \cite{MKM99}.}
\end{figure}

Figure \ref{Figure1a} shows the normalized mean streamwise velocity; 
note the almost 
perfect overlapping of the results corresponding to the models tested. 
We interpret this behavior as a measure of our success in enforcing a 
constant mass flux through the channel.

\begin{figure}
\centerline{
\hbox{
\psfig{figure=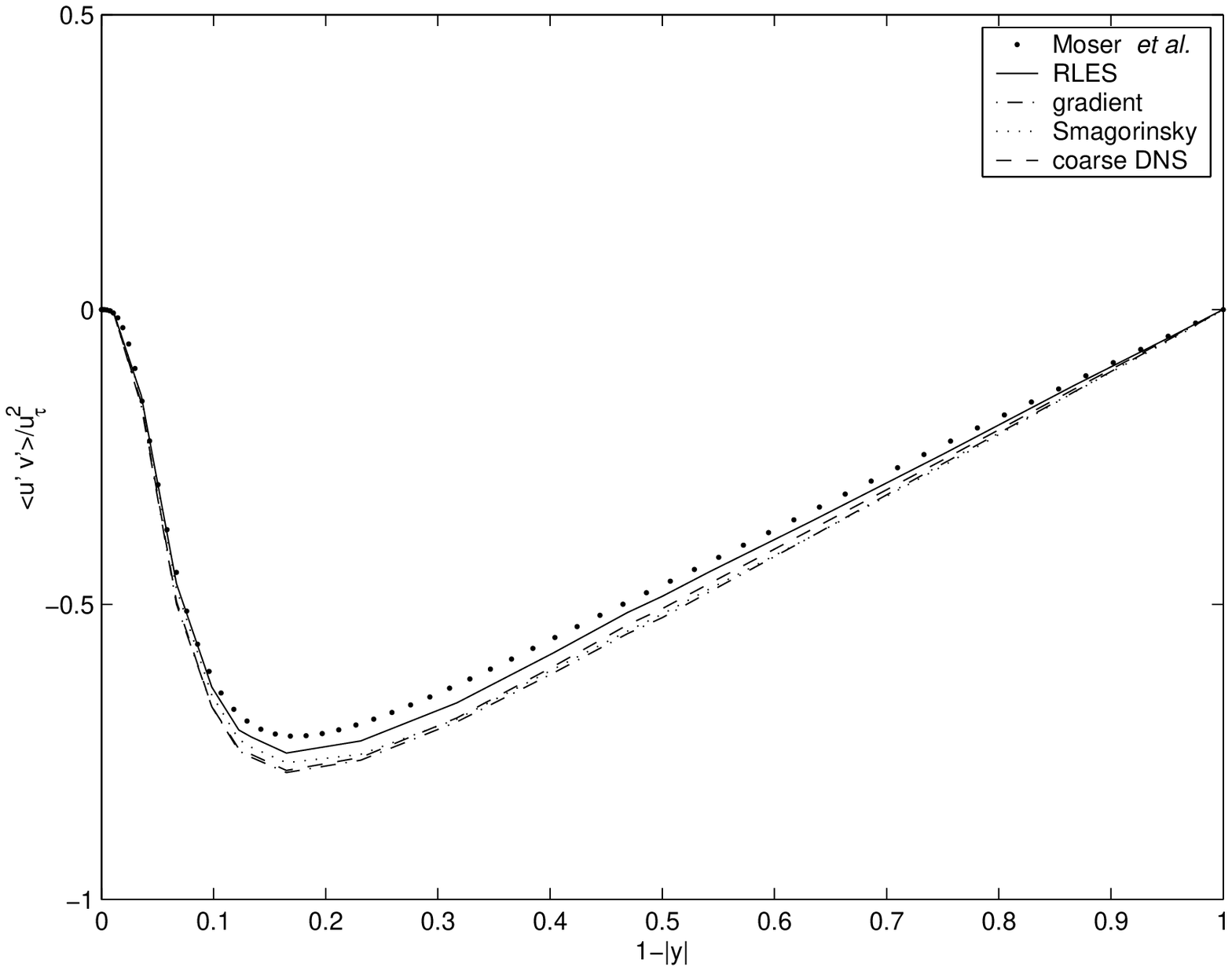,height=300pt,width=400pt,angle=0}}
}
\caption{\label{Figure2a} The $x,y$-component of the Reynolds stress, $Re_{\tau}=180$.
We compared the RLES model (\ref{rational.LES}), the gradient model
(\ref{gradient.model}), the Smagorinsky model, and a coarse DNS,
with the fine DNS of Moser, Kim, and Mansour \cite{MKM99}.}
\end{figure}

Figure \ref{Figure2a} presents the normalized $x,y$-component of the 
Reynolds stress.
The RLES model (\ref{rational.LES}) is a clear improvement over the rest 
(i.e., closest to the fine DNS).

\begin{figure}
\centerline{
\hbox{
\psfig{figure=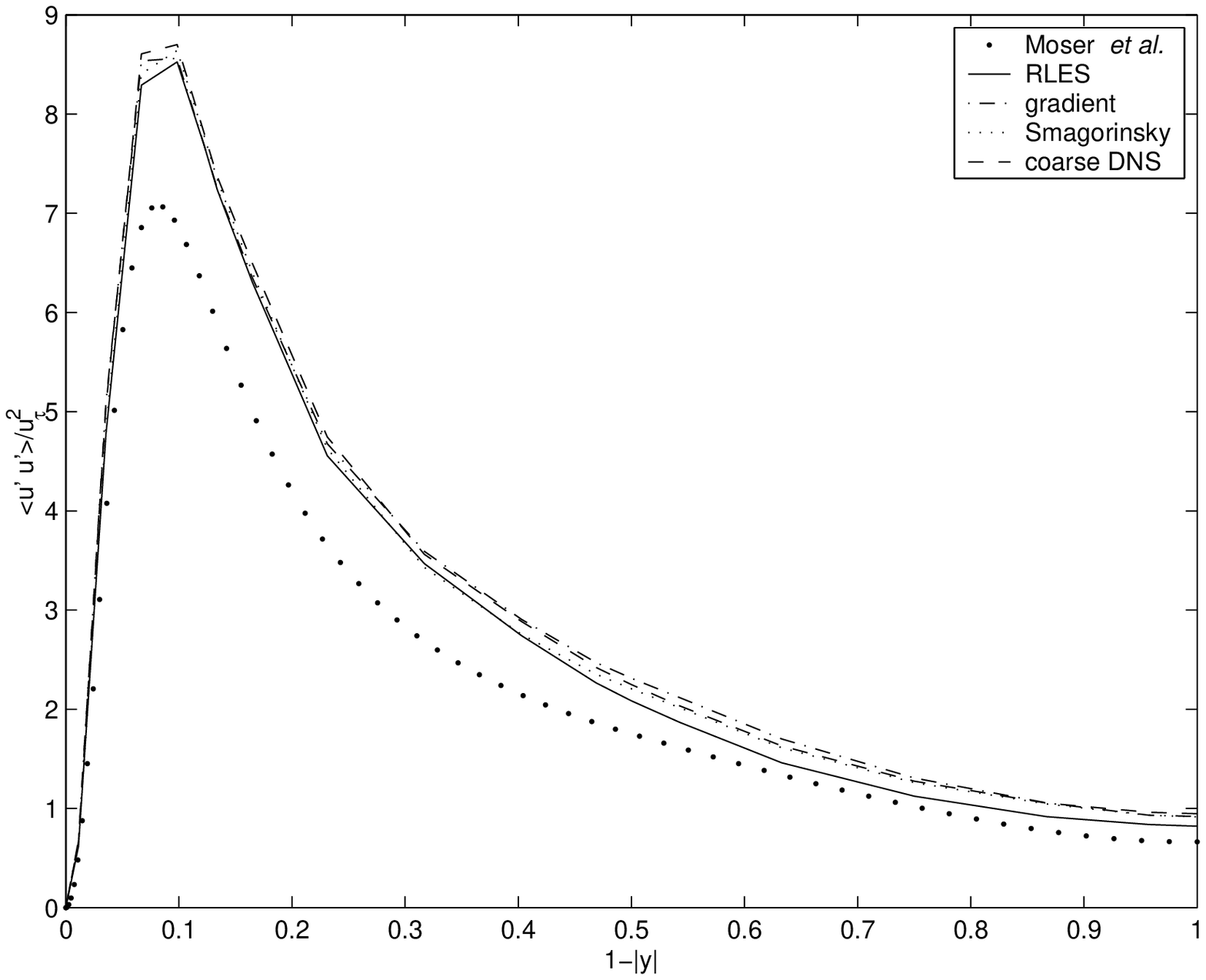,height=300pt,width=400pt,angle=0}}
}
\caption{\label{Figure3a} Rms values of streamwise velocity fluctuations, $Re_{\tau}=180$.
We compared the RLES model (\ref{rational.LES}), the gradient model
(\ref{gradient.model}), the Smagorinsky model, and a coarse DNS,
with the fine DNS of Moser, Kim, and Mansour \cite{MKM99}.}
\end{figure}

Similarly, the normalized rms values of the streamwise velocity 
fluctuations in Figure \ref{Figure3a} 
show a better (closer to the fine DNS benchmark results in \cite{MKM99}) 
behavior for the RLES model (\ref{rational.LES}).

\begin{figure}
\centerline{
\hbox{
\psfig{figure=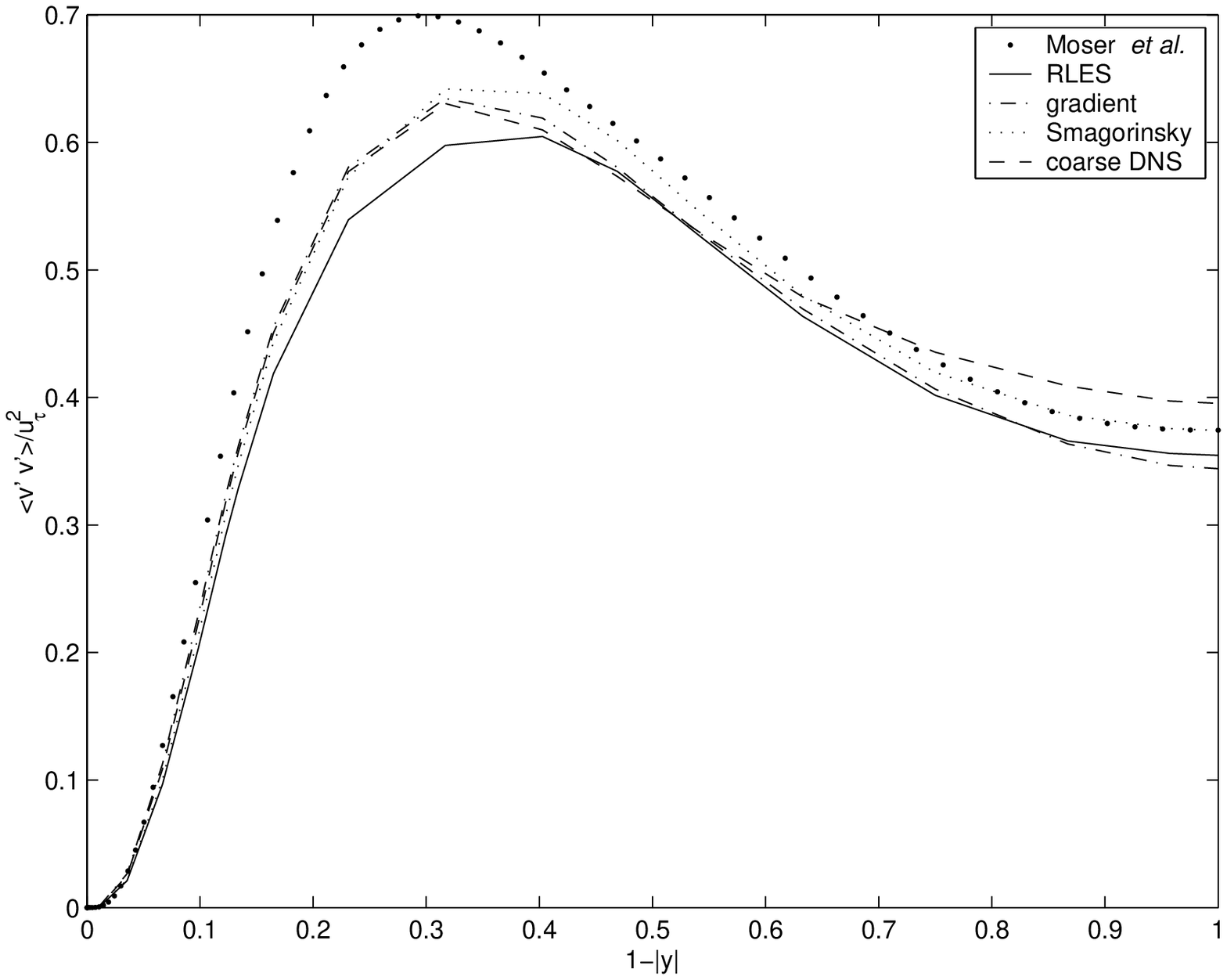,height=300pt,width=400pt,angle=0}}
}
\caption{\label{Figure4a} Rms values of wall-normal velocity fluctuations, $Re_{\tau}=180$.
We compared the RLES model (\ref{rational.LES}), the gradient model
(\ref{gradient.model}), the Smagorinsky model, and a coarse DNS,
with the fine DNS of Moser, Kim, and Mansour \cite{MKM99}.}
\end{figure}

\begin{figure}
\centerline{
\hbox{
\psfig{figure=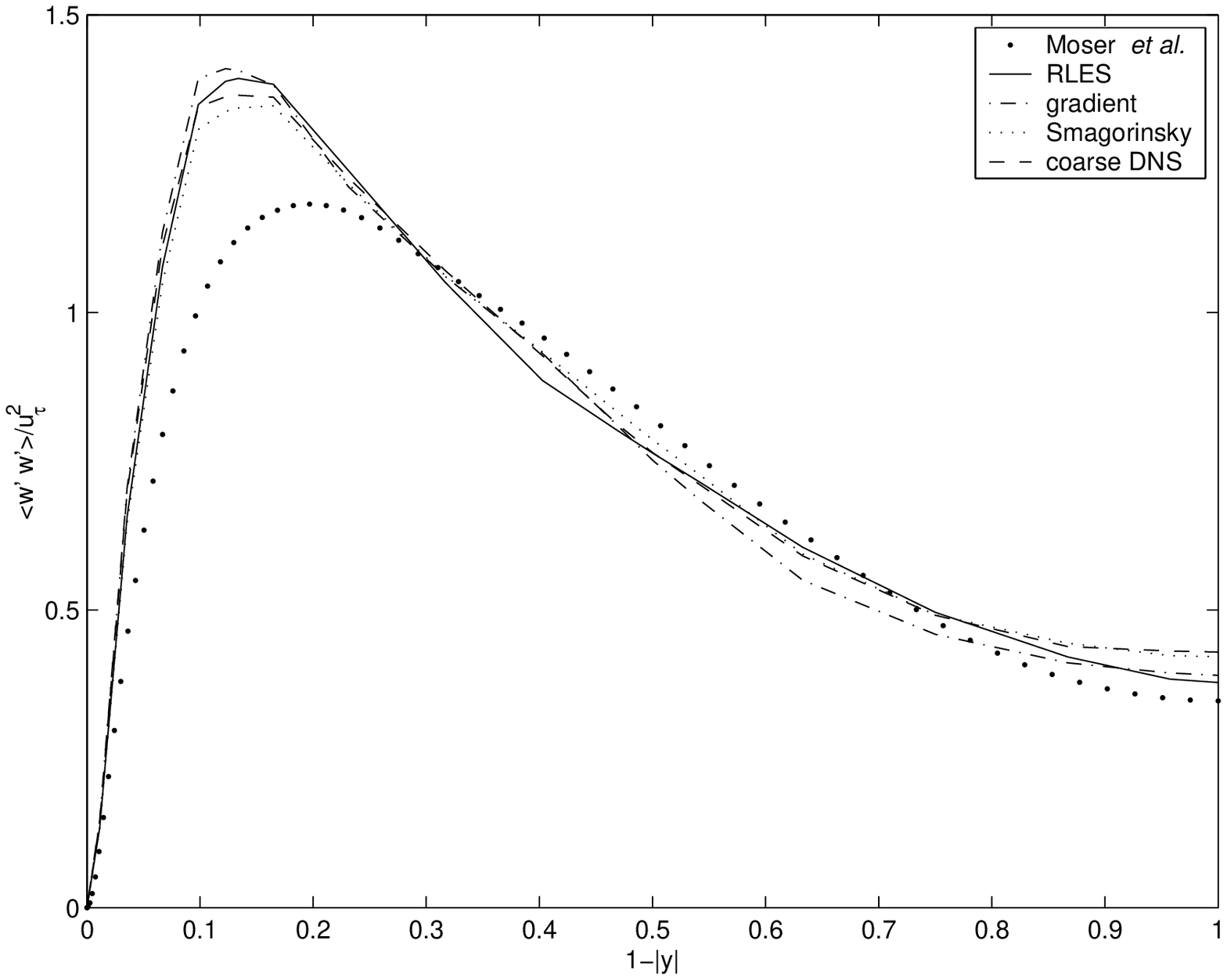,height=300pt,width=400pt,angle=0}}
}
\caption{\label{Figure5a} Rms values of spanwise velocity fluctuations, $Re_{\tau}=180$.
We compared the RLES model (\ref{rational.LES}), the gradient model
(\ref{gradient.model}), the Smagorinsky model, and a coarse DNS,
with the fine DNS of Moser, Kim, and Mansour \cite{MKM99}.}
\end{figure}

Figures \ref{Figure4a} and \ref{Figure5a}, containing the 
normalized rms values for the wall-normal and spanwise 
velocity fluctuations, merit a more detailed discussion.   
Figure \ref{Figure4a} shows the normalized rms values of the 
wall-normal velocity 
fluctuations.  Here, the RLES model (\ref{rational.LES})
performs worse than the gradient model (\ref{gradient.model}) near the wall.
Away from the wall, the gradient and the RLES model perform similarly.
The best results are obtained with the Smagorinsky model.

The normalized rms values of the spanwise velocity fluctuations 
in Figure \ref{Figure5a} are better for the RLES model (\ref{rational.LES})
than for the gradient model, except for a portion of $1-|y|$
roughly between 0.3 and 0.5. 
The Smagorinsky model gives the best results near the wall, but it
overpredicts the correct value near the center of the channel.

In conclusion, the RLES model (\ref{rational.LES}) performs better than 
the gradient model, with the exception of the normalized rms values 
of the wall-normal velocity fluctuations.
The RLES model (\ref{rational.LES}) is also more stable numerically 
than the gradient model.

\section{{\it A Posteriori} Tests for $Re_{\tau}=395$}

\begin{figure}
\centerline{
\hbox{
\psfig{figure=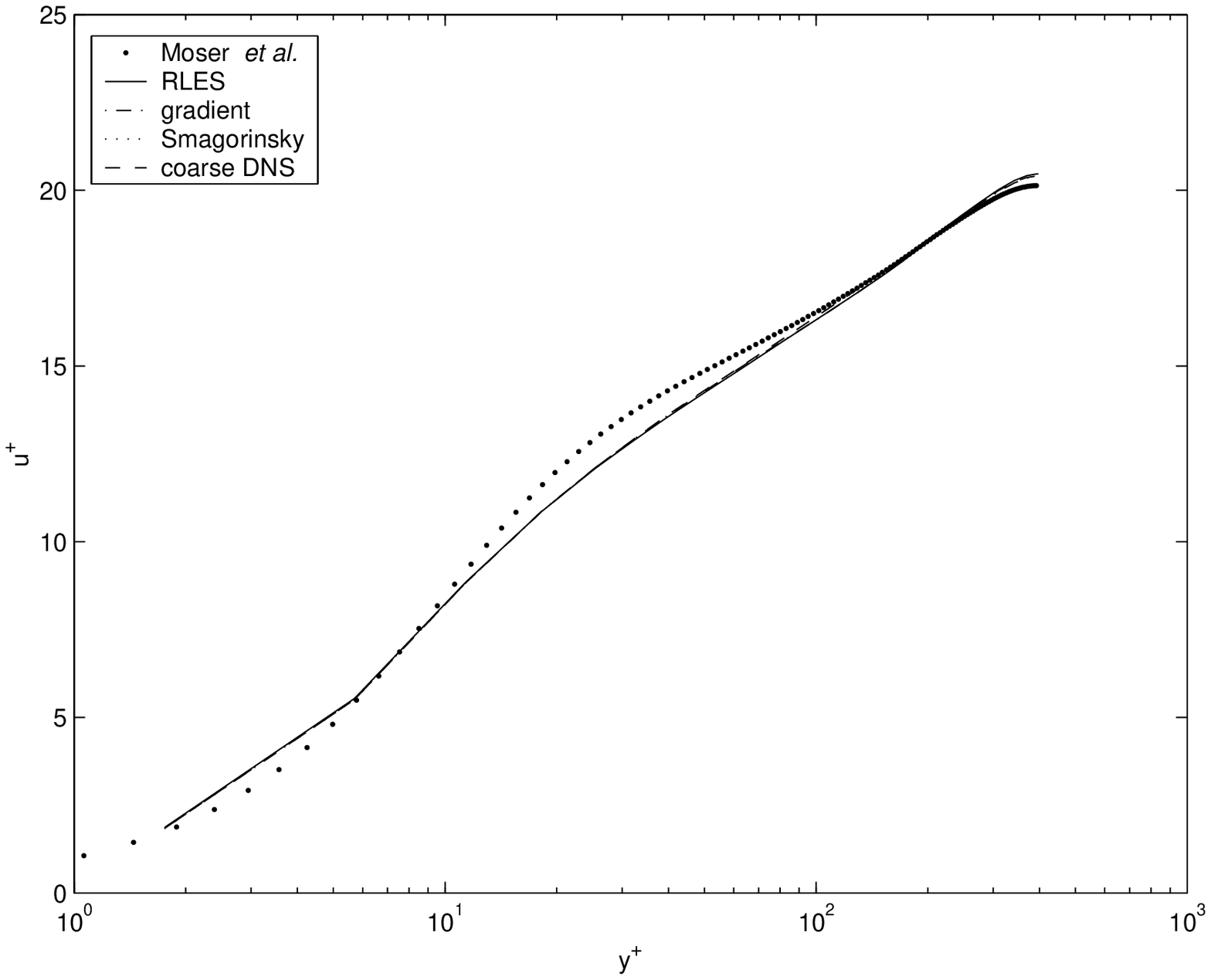,height=300pt,width=400pt,angle=0}}
}
\caption{\label{Figure1b} Mean streamwise velocity, $Re_{\tau}=395$.
We compared the RLES model (\ref{rational.LES}), the gradient model
(\ref{gradient.model}), the Smagorinsky model, and a coarse DNS,
with the fine DNS of Moser, Kim, and Mansour \cite{MKM99}.}
\end{figure}

We ran simulations with all four models for  $Re_{\tau}=395$, and 
we compared our results with the fine DNS in \cite{MKM99}.
Again, as in the $Re_{\tau}=180$ case, the normalized mean
streamwise velocity fluctuations in Figure \ref{Figure1b}
are practically identical; this time, however, they do
not overlap onto that corresponding to the fine DNS.
Nevertheless, the mean flows are the same, and this is 
supported by the fact that the models underpredict the 
correct value near the wall but overpredict it
away from the wall.

\begin{figure}
\centerline{
\hbox{
\psfig{figure=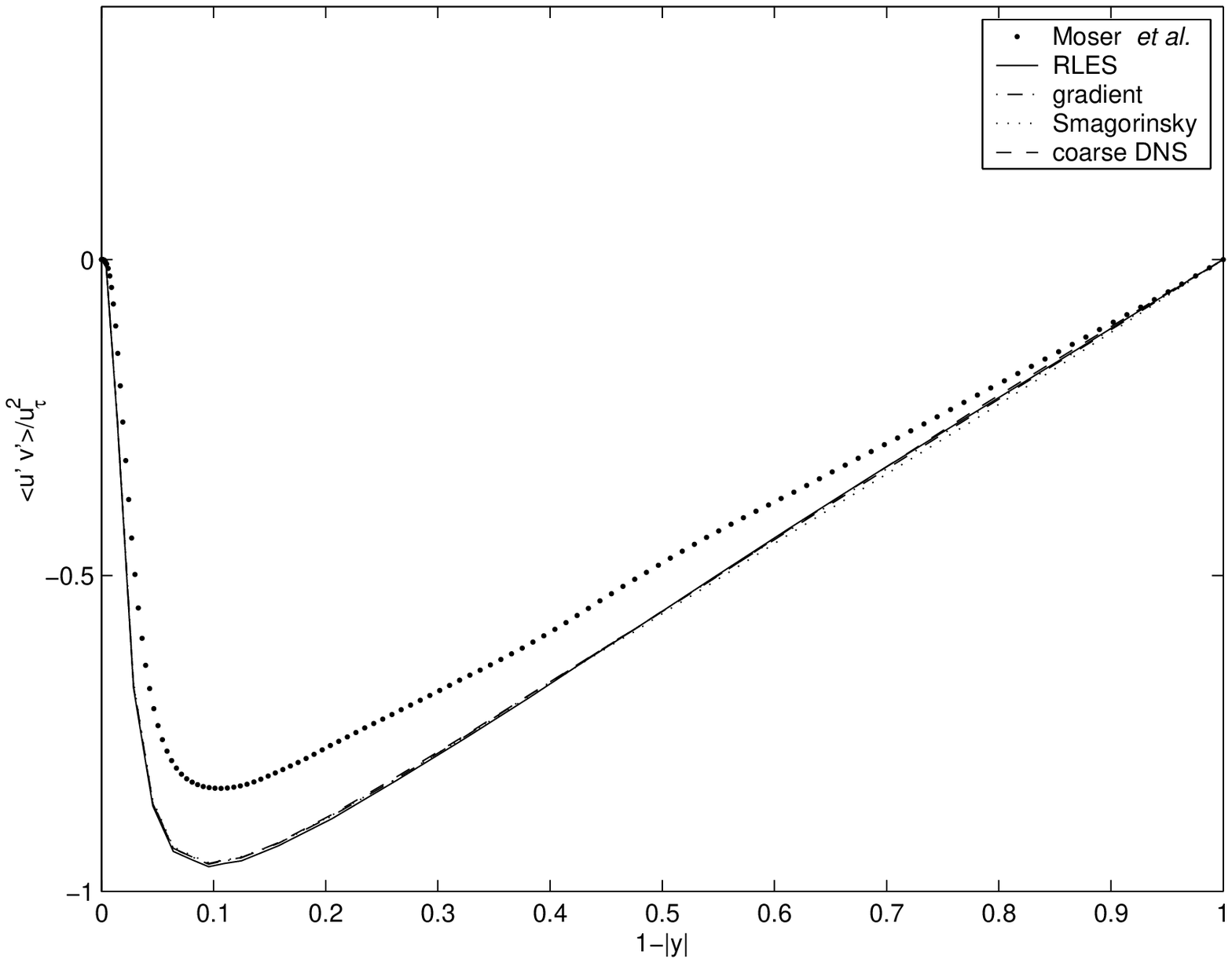,height=300pt,width=400pt,angle=0}}
}
\caption{\label{Figure2b} The $x,y$-component of the Reynolds stress, $Re_{\tau}=395$.
We compared the RLES model (\ref{rational.LES}), the gradient model
(\ref{gradient.model}), the Smagorinsky model, and a coarse DNS,
with the fine DNS of Moser, Kim, and Mansour \cite{MKM99}.}
\end{figure}

The normalized $x,y$-component of the Reynolds stress
in Figure \ref{Figure2b} is almost identical for all four models.
This behavior was also noticed by Winckelmans et al.
\cite{WWVJ01}.

\begin{figure}
\centerline{
\hbox{
\psfig{figure=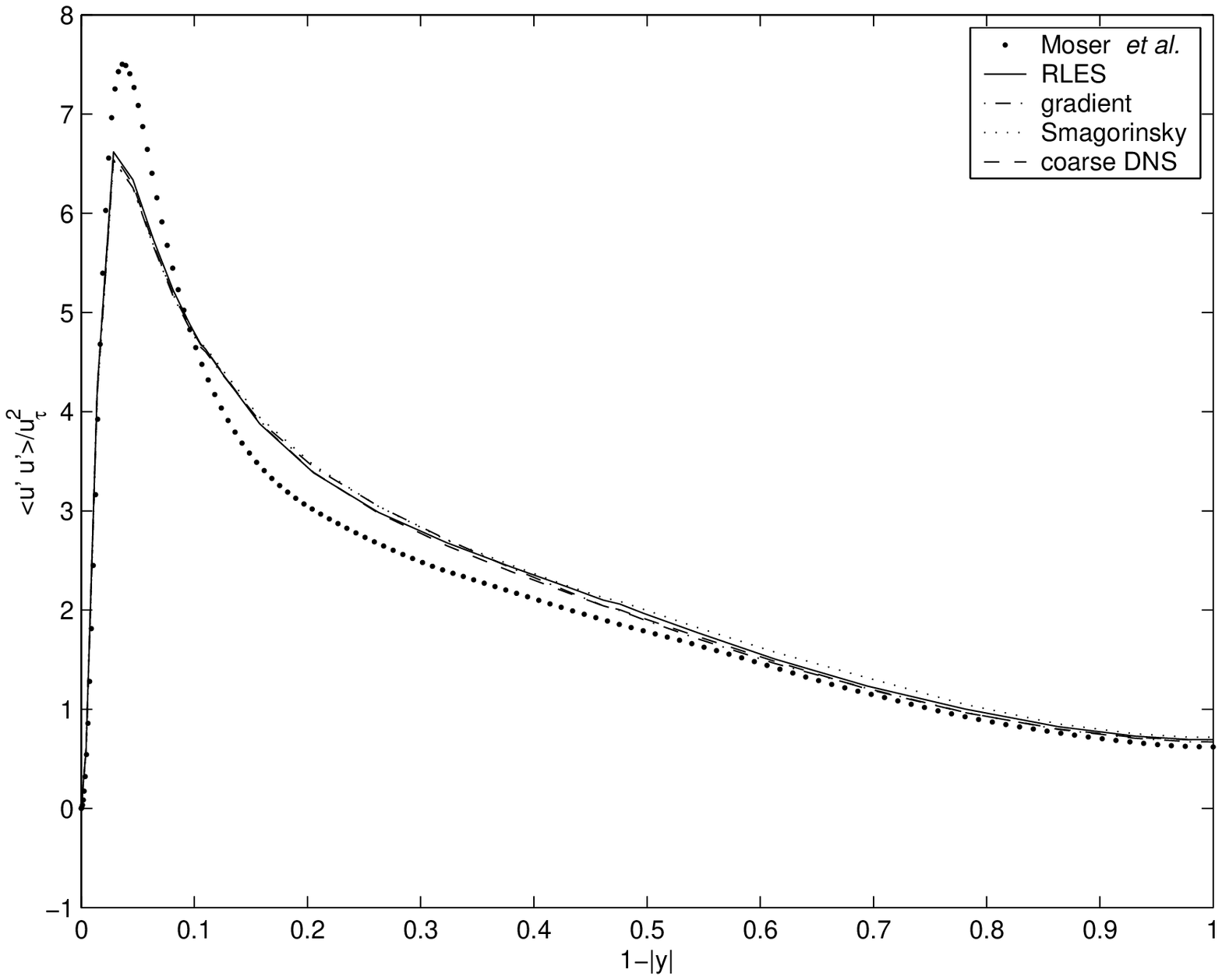,height=300pt,width=400pt,angle=0}}
}
\caption{\label{Figure3b} Rms values of streamwise velocity fluctuations, $Re_{\tau}=395$.
We compared the RLES model (\ref{rational.LES}), the gradient model
(\ref{gradient.model}), the Smagorinsky model, and a coarse DNS,
with the fine DNS of Moser, Kim, and Mansour \cite{MKM99}.}
\end{figure}

The same behavior can be noticed for the normalized rms values
of the streamwise velocity fluctuations in Figure \ref{Figure3b}:
the profiles for the four models are almost identical.
The RLES model (\ref{rational.LES}) performs slightly  better 
near the wall, and the gradient model performs slightly better 
away from the wall.
They both perform better than the Smagorinsky model near the 
center of the channel.

\begin{figure}
\centerline{
\hbox{
\psfig{figure=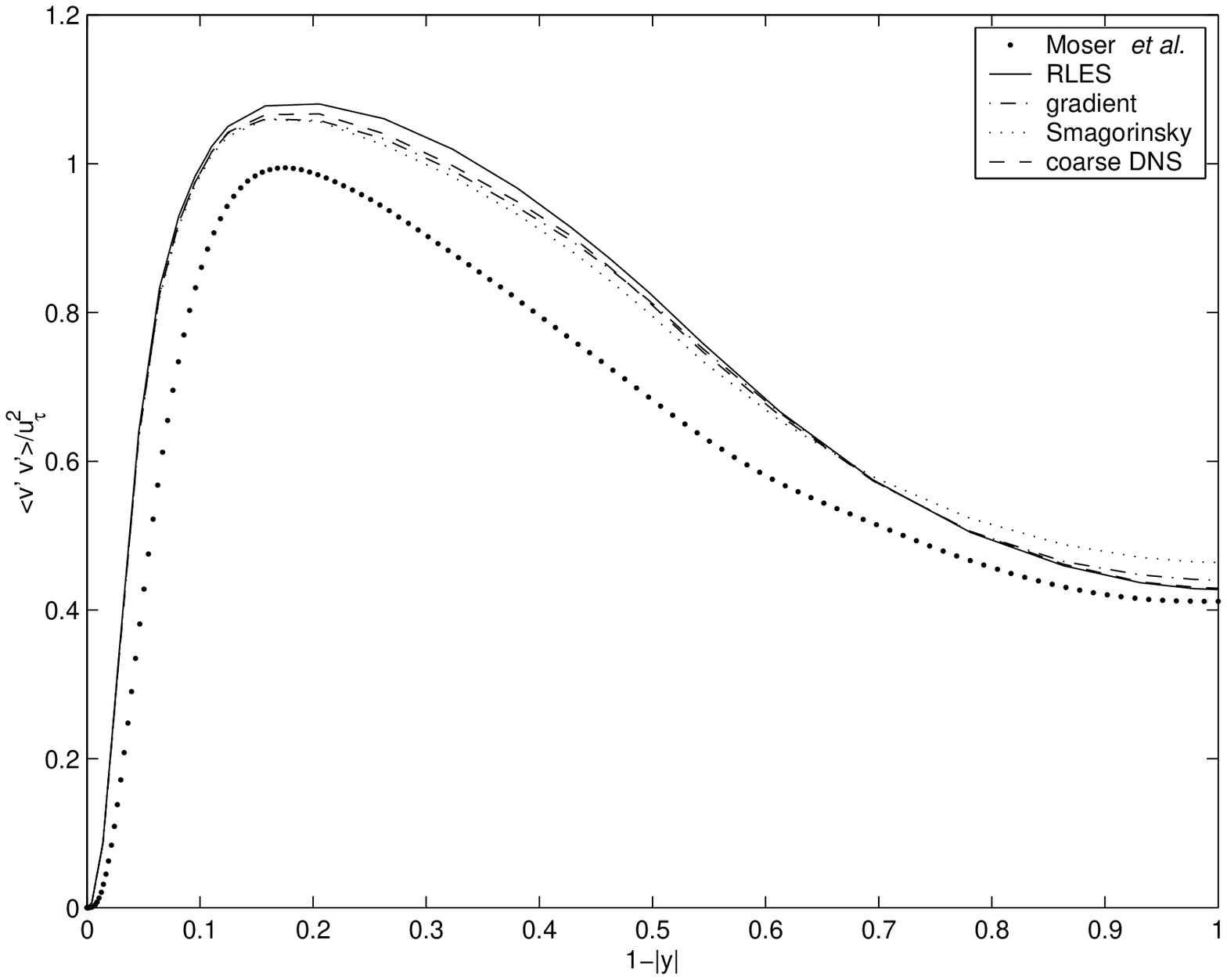,height=300pt,width=400pt,angle=0}}
}
\caption{\label{Figure4b} Rms values of wall-normal velocity fluctuations, $Re_{\tau}=395$.
We compared the RLES model (\ref{rational.LES}), the gradient model
(\ref{gradient.model}), the Smagorinsky model, and a coarse DNS,
with the fine DNS of Moser, Kim, and Mansour \cite{MKM99}.}
\end{figure}

Figure \ref{Figure4b} presents the normalized rms values of the 
wall-normal velocity fluctuations.
As for the $Re_{\tau}=180$ case, the gradient model performs 
better near the wall, and the RLES model (\ref{rational.LES}) 
performs better away from the wall.
The Smagorinsky model performs best, but it overpredicts the 
correct value near the center of the channel.

\begin{figure}
\centerline{
\hbox{
\psfig{figure=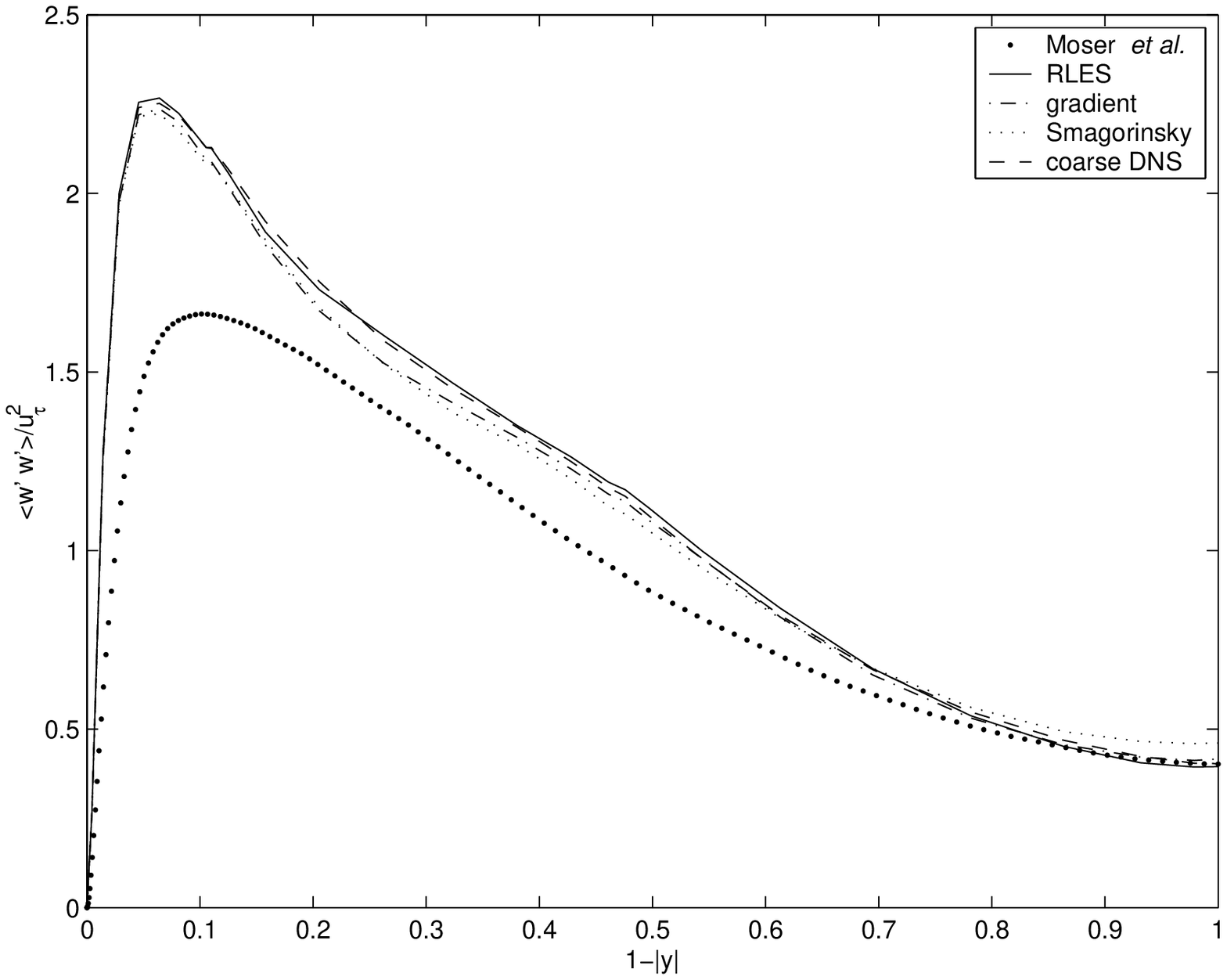,height=300pt,width=400pt,angle=0}}
}
\caption{\label{Figure5b} Rms values of spanwise velocity fluctuations, $Re_{\tau}=395$.
We compared the RLES model (\ref{rational.LES}), the gradient model
(\ref{gradient.model}), the Smagorinsky model, and a coarse DNS,
with the fine DNS of Moser, Kim, and Mansour \cite{MKM99}.}
\end{figure}

The same behavior is observed for the  normalized rms values
of the spanwise velocity fluctuations in Figure \ref{Figure5b},
and again the Smagorinsky model performs best.

In conclusion, for the $Re_{\tau}=395$ case, the gradient and 
the RLES model (\ref{rational.LES}) yield comparable results.
The best results, however, are obtained by using the 
Smagorinsky model.

Again, as in the $Re_{\tau}=180$ case, the RLES model 
(\ref{rational.LES}) is much more stable numerically than the 
gradient model.

\section{Conclusions} 
We have used a spectral element code to test the 
RLES model (\ref{rational.LES})  
in the numerical
simulation of incompressible channel flows at 
$Re_{\tau}=180$ and $Re_{\tau}=395$.
This approximate deconvolution model is
based on a rational (Pad\'e) approximation to the Fourier 
transform of the Gaussian filter and is proposed as
an alternative to the gradient model 
(\ref{gradient.model}).
We compared the RLES model (\ref{rational.LES}) with 
the gradient model, the Smagorinsky model, and a 
coarse DNS (no LES model).
The corresponding results were benchmarked against 
the fine DNS calculation of Moser, Kim, and Mansour 
\cite{MKM99}.

The RLES model (\ref{rational.LES}) yielded the best
results for the $Re_{\tau}=180$ case.
These improved results were accompanied by a much
increased numerical stability compared with the 
gradient model.

The situation was different for the  $Re_{\tau}=395$ case.
Here the RLES model (\ref{rational.LES}) and the gradient 
model yielded comparable results, and the Smagorinsky
model performed the best.
Again, the RLES model (\ref{rational.LES}) was much
more stable numerically than the gradient model.

We believe that these results for the RLES 
model are encouraging.
They also support our initial thoughts:
The RLES model is an improvement over the gradient
model as a subfilter-scale model, and this is illustrated 
by the improved results for the $Re_{\tau}=180$ case.
The RLES model is also more stable numerically because of the
additional smoothing operator, and this feature is manifest 
for both low ($Re_{\tau}=180$) and moderate ($Re_{\tau}=395$)
Reynolds number flows.

However, the RLES model accounts just for the subfilter-scale
part of the stress reconstruction.
The information lost at the subgrid-scale level must be 
accounted for in a different way, as advocated by 
Carati et al. \cite{CWJ01}.
This was illustrated by the fact that, even for a moderate 
Reynolds number ($Re_{\tau}=395$) flow, the Smagorinsky
model, a classical eddy-viscosity model, performed best.

Along these lines, our next step will be to develop a 
mixed model, consisting of the RLES model supplemented
by a Smagorinsky model.
We also plan to study improved boundary conditions, 
the commutation error \cite{GM95}, \cite{Gho96}, 
and the relationship between the filter radius and 
the mesh-size in a spectral element discretization.

\bigskip
\bigskip
\bigskip

{\large {\bf Acknowledgments.}} This work was supported in part by the Mathematical, 
               Information, and Computational Sciences Division subprogram of the 
               Office of Advanced Scientific Computing  Research, U.S. Dept. of 
               Energy, under Contract W-31-109-Eng-38.

We thank Professor R. Moser and Mr. A. Das
for helpful communications that improved this paper.

\newpage

\vspace*{8.0cm}

\newpage

\vspace*{5.0cm}  

\begin{table}[h]
\small
\centering
\caption{Parameters for the numerical simulations.}
\vspace*{0.5cm}  
\label{table.param}
\begin{tabular}{c c c}
\hline
\hline
  Nominal $Re_{\tau}$  &  $L_x \times L_y \times L_z$     &  $N_x \times N_y \times N_z$ \\ 
\hline
\ \\
  180  & $4 \pi \times 2 \times \frac{4}{3} \pi$  &  $36 \times 37 \times 36$ \\
  395  & $2 \pi \times 2 \times \pi$              &  $72 \times 55 \times 54$ \\
\ \\
\hline
\hline
\end{tabular} 
\end{table}

\newpage

\vspace*{3.0in}

\newpage

\vspace*{3.0in}

\newpage

\vspace*{2.0in}

\newpage

\vspace*{2.0in}

\newpage

\vspace*{5.0cm}  

\begin{table}[h]
\small
\centering
\caption{Computed $u_\tau$ and $Re_\tau$.}
\vspace*{0.5cm}  
\label{table.utau}
\begin{tabular}{c c c c c}
\hline
\hline
  Fixed $U_m$   &   Nominal $Re_\tau$  &  Case  &  Computed $u_\tau$  &  Computed $Re_\tau$\\ 
\hline
\ \\
  15.63   &   180   &   RLES         &  0.9879448  &  177.8352\\
          &         &   gradient     &  0.9890118  &  178.0222\\
          &         &   Smagorinsky  &  0.9917144  &  178.5120\\
          &         &   coarse DNS   &  0.9873800  &  177.7184\\
\ \\
\hline
\ \\
  17.54   &   395   &   RLES         &  1.001025319   &  395.4071960\\
          &         &   gradient     &  1.005021334   &  396.9859924\\
          &         &   Smagorinsky  &  0.9974176884  &  393.9718933\\
          &         &   coarse DNS   &  0.9901855588  &  391.1294861\\
\ \\
\hline
\hline
\end{tabular} 
\end{table}

\newpage

\vspace*{2.0in}

\newpage

\vspace*{2.0in}

\newpage

\vspace*{2.0in}

\newpage

\vspace*{2.0in}

\newpage

\vspace*{2.0in}

\newpage

\vspace*{2.0in}

\newpage

\vspace*{2.0in}

\newpage

\vspace*{2.0in}

\newpage

\vspace*{2.0in}

\newpage

\vspace*{2.0in}

\end{document}